\documentclass[12pt]{amsart}
\usepackage{amsmath,amssymb}
\usepackage{enumerate}
\def\B{\Bbb B}
\def\C{\Bbb C}
\def\D{\Bbb D}
\def\Om{\Omega}
\def\CC{\mathcal C}
\def\var{\varepsilon}

\newtheorem{thm}{Theorem}
\newtheorem{corollary}[thm]{Corollary}

\def\ds{\displaystyle}

\begin{document}

\title[Estimates for the squeezing function with applications]
{Estimates for the squeezing function near strictly pseudoconvex boundary points with applications}

\author{Nikolai Nikolov}
\address{Institute of Mathematics and Informatics\\Bulgarian Academy
of Sciences\\ Acad. G. Bonchev 8, 1113 Sofia, Bulgaria
\vspace{0.1cm}
\newline Faculty of Information Sciences\\
State University of Library Studies and Information Technologies\\
Shipchenski prohod 69A, 1574 Sofia,
Bulgaria}

\email{nik@math.bas.bg}

\author{Maria Trybu\l a}

\address{Faculty of Mathematics and Informatics\\Adam Mickiewicz University\\ Umultowska 87, 61-614 Pozna\'{n}, Poland}

\email{maria.h.trybula@gmail.com}

\thanks{The first named author was partially supported by the Bulgarian National Science Fund,
Ministry of Education and Science of Bulgaria under contract DN 12/2.
The second named author was supported by National Center of Science (Poland), grant no. 2013/10/A/ST1/00091.
Part of this work was done while her stay in Sofia during {\it NTADES 2018.}}



\begin{abstract}
An extension of the estimates for the squeezing function of strictly pseudoconvex domains
obtained recently by J. E. Forn\ae ss and E. Wold in \cite{FW1} is applied to
derive sharp boundary behaviours of invariant metrics and Bergman curvatures.
\end{abstract}

\maketitle

\noindent{\bf I.} Let $\Om$ be a bounded domain in $\C^n.$ Following \cite{DGZ1} (which is motivated by \cite{LSZ1, Yeu}),
we define the squeezing function $\sigma_\Om$ as follows: For any injective holomorphic map $f:\Om\rightarrow\B_n$ with $f(z)=0$
set $$\sigma_{\Om,f}(z)=\max\{r>0:r\B_n\subset f(\Om)\},$$
where $\B_n$ is the unit ball in $\C^n,$ and then
$$\sigma_\Om(z)=\sup_f\sigma_{\Om,f}(z).$$

Roughly speaking, $\sigma_\Om(z)$ describes how much $\Om$ looks like the unit ball at the point $z.$ Though the definition of the squeezing function is simple, it turns out that many geometric and analytic properties of $\Om$ are hidden in $\sigma_\Om.$ Some basic properties of $\sigma_\Om$ are gathered in \cite{DGZ1,DGZ2}, for example, $\sigma_\Om$ is continuous and $\sup_f$ is always attained; $\sigma_\Om(z)=1$ if and only if $\Om$ is biholomorphic to $\B_n.$

It was proven in \cite{DGZ2} (see also \cite{KZ}) that
$$\lim_{z\rightarrow\partial\Om}\sigma_\Om(z)=1$$
if $\Om$ is a $\CC^2$-smooth strictly pseudoconvex domain. This result was improved in \cite{FW1} when the boundary has higher regularity. More precisely, it was shown that if $b\Om$ is $\CC^3$-smooth (resp. $\CC^4$-smooth), then
$$\sigma_\Om(z)\geq 1-C\sqrt{\delta_{\Om}(z)}\ (\mbox{resp.}\ \sigma_\Om(z)\geq 1-C\delta_{\Om}(z)),$$
where $C>0$ is a constant and $\delta_\Om(z)=d(z,\partial\Omega).$

E.F. Wold posed the question what is the optimal estimate for $\sigma_\Om$ if $b\Om$ is $\CC^{2,\var}$-smooth, see \cite{Wold} (see also \cite{FW2}).
Theorem \ref{main} below shows that similar estimates as above hold in the $\CC^{2,\var}$- and $\CC^{3,\var}$-smooth cases.

Consider the following assumption:

\vspace{0.2cm}
$(\ast)$ $\var\in (0,1],$ $k\in\{0,1\},$ $\var_k=\frac{k+\var}{2},$ and $p$ is a $\CC^{k+2,\var}$-smooth
strict pseudoconvex boundary point of a domain $\Om$ in $\C^n.$

\vspace{0.1cm}
\begin{thm}\label{main}
If $(\ast)$ holds and $\overline{\Om}\Subset \C^n$ admits a Stein neighborhood basis, then
there exist a constant $C>0$ and a neighborhood $U$ of $p$ such that
$$\sigma_\Om(z)\geq 1-C\delta_\Om(z)^{\var_k},\quad z\in\Om\cap U.$$
\end{thm}

Note that no global strict pseudoconvexity is assumed.

\begin{proof} The proof follows the same lines as that of \cite[Theorem 1.1]{FW1} (the $\CC^3$- and $\CC^4$-smooth
cases), replacing there $k$ by $2+2\var_k$ and setting $\tilde\eta=\frac{\eta}{1-C\eta^{\frac{\var_k}{1-\var_k}}},$
$\eta=(1-r)^{1-\var_k}$ ($0<\tilde\eta=\eta\ll 1$ is fixed if $\var_k=1$).
\end{proof}

Similarly to \cite{DF,FW1,Zha}, we will apply Theorem \ref{main} to derive sharp boundary behaviours of invariant metrics and Bergman curvatures.
\smallskip

\noindent{\bf II.} D. Ma in \cite{Ma} obtained sharp estimates for the behaviour of the Kobaya\-shi metric near $\CC^3$-smooth strictly pseudoconvex points.
S. Fu in \cite{Fu} refined these estimates in the $\CC^\infty$-smooth case. It was shown in \cite{Nik} that Ma's approach works also in the $\CC^{2,\var}$-smooth case. Similarly, it can be observed that Fu's arguments can be adapted to the $\CC^{3,\var}$-smooth case. Before we state the precise result, recall the definitions of the Carath\'eodory and the Kobayashi metrics:
$$\gamma_{\Om}(z;X)=\sup\{|f'(z)X|:f\in\mathcal{O}(\Om, \D)\},$$
$$\kappa_{\Om}(z;X)=\inf\{|\alpha|:\exists\varphi\in \mathcal{O}(\D,\Om),\,\varphi(0)=z,\,\alpha\varphi'(0)=X\},$$
where $\D$ is the unit disc.

\begin{thm}\label{nik} (D. Ma $\&$ S. Fu $+$) If $(\ast)$ holds, then there exist a constant $C>0$ and a neighborhood $U$ of $p$ such that
\begin{multline*}
(1-C\delta_\Om(z)^{\var_k})\Bigl(\frac{|\langle\partial\delta_\Om(z),X\rangle|^2}{\delta^2_\Om(z)}+\frac{L(z;X)}{\delta_\Om(z)}\Bigr)^{1/2}
\leq\kappa_\Om(z;X) \\
\leq (1+C\delta_\Om(z)^{\var_k})\Bigl(\frac{|\langle\partial\delta_\Om(z),X\rangle|^2}{\delta^2_\Om(z)}+\frac{L(z;X)}{\delta_\Om(z)}\Bigr)^{1/2},
\quad z\in\Om\cap U,\ X\in\C^n,
\end{multline*}
where $L$ is the Levi form of $-\delta_\Om.$
\end{thm}

Theorem \ref{nik} remains true in the $\CC^{2}$-smooth case, replacing the term
$C\delta_\Om(z)^{\var_k}$ by any positive number.

Setting $\varepsilon=1,\,k=1,$ we rediscover \cite[Theorem 2.2]{Fu}, which is more precise in some sense
than \cite[Theorem B]{Ma} despite of \cite[Example on p. 337]{Ma}. This is due to the fact that we treat
the holomorphic tangent vectors as vectors of the interior points rather than projecting them to the boundary,
decomposing into complex normal and tangential components.

\begin{proof}[Sketch of the proof of Theorem \ref{nik}]
After translations and rotations, we may assume that $0$ is the nearest point to $p=$ on $\partial\Om$ and the negative $\textup{Im}\,z_1$-axis is the outside normal direction at $0.$

First, we approximate $\partial\Om$ near the origin by a biholomorphic image of a ball to the ($2+2\varepsilon_k$)-th order, that is $\partial\Omega$ near $0$ takes the following form: $
\textup{Re}\,\zeta_1=|\zeta'|^2+O(|\zeta'|^{2+2\varepsilon_k}+|\textup{Im}\ \zeta_1|^{2+2\varepsilon_k}).$  This is just a weak version of \cite[Lemma 2.1]{Fu}. For that it suffices to use the map $\Phi_1$ in the $\mathcal{C}^{2,\varepsilon}$-smooth case (resp. the maps $\Phi_1$ and $\Phi_2$ in the $\mathcal{C}^{3,\varepsilon}$-smooth case).

Second, we follow the proof of \cite[Theorem 2.2]{Fu} with one change resulting from the lower regularity: we replace the exponent $2$ in the inequality $|M(\zeta)|\leq C|\zeta_1|^2$ in the middle of \cite[p. 203]{Fu} by $\varepsilon_k,$ and then the inequality $\big{|}\lVert\Psi(\zeta)\rVert^2-1\big{|}<C\tau$ becomes $\big{|}\lVert\Psi(\zeta)\rVert^2-1\big{|}\leq C\tau^{\varepsilon_k-1}.$
\end{proof}

The next corollary is a consequence of Theorems \ref{main} and \ref{nik}, and the following estimates
(see \cite{DGZ1,FW1,NA}):
$$\sigma_\Om\le\frac{\gamma_\Om}{\kappa_\Om}\le 1,\quad\sigma_\Om^{n+1}\le\frac{1}{\sqrt{n+1}}
\cdot\frac{\beta_\Om}{\kappa_\Om}\le\sigma_\Om^{-n-1}.$$

\begin{corollary}\label{cor} If $(\ast)$ holds and $\overline{\Om}\Subset\C^n$ admits a Stein neighborhood basis,
then the estimates from Theorem \ref{nik} remain true for $\gamma_\Om$ and $\frac{\beta_\Om}{\sqrt{n+1}}$ instead
of $\kappa_{\Om},$ where $\beta_\Om$ is the Bergman metric \textup{(}see below\textup{)}.
\end{corollary}

When $\Om$ is a $\CC^{3}$- or $\CC^{4}$-smooth strictly pseudoconvex domain,
Corollary \ref{cor} can be found in \cite{FW1,DF}. We also refer to \cite{Nik}
in the case of $\beta_\Om$ and $\CC^{2,\var}$-smoothness. There the global assumption
is weaker, namely pseudoconvexity. Later, we will obtain the same in the $\CC^{3,\var}$-smooth
case, $\var\neq 1.$

Note also that when $\Om$ is a $\CC^\infty$-smooth strictly pseudoconvex domain,
the estimate from Theorem \ref{main} implies the same estimate for $\gamma_\Om.$ This fact is proved in
\cite{Fu} by using the Fornaess embedding theorem (see e.g. \cite{JP}), where the global strictly pseudoconvexity
is essential.
\smallskip

\noindent{\bf III.} As another application  of Theorem \ref{main}, we will explore the boundary behavior of holomorphic
sectional curvature $s_\Om$, the Ricci curvature $r_\Om,$ and the scalar curvature $t_\Om$ of the Bergman metric
$\beta_\Om$ near $\CC^{2,\var}$- and $\CC^{3,\var}$-smooth strictly pseudoconvex boundary points of an arbitrary
pseudoconvex domain $\Om$ in $\C^n$ (even not necessarily bounded).
When $\Om$ is $\mathcal C^\infty$-smooth and strictly pseudoconvex, more precise results can be obtained by using
Fefferman's asymptotic expansion of the Bergman kernel (see \cite{Eng}).

Our approach is based on the fact that these Bergman invariants can be expressed by extremal domain functions.

Denote by $L^2_h(\Om)$ all holomorphic square integrable functions on $\Om;$ we write $||\cdot||_\Om$
for the norm in $L^2(\Om).$ Set:
$$J_\Om^0(z)=\sup\{|f(z)|^2:f\in L^2_h(\Om),||f||_D\leq 1\},$$
$$J_\Om^1(z;X)=\sup\{|f'(z)X|^2:f\in L^2_h(\Om),||f||_D\leq 1,f(z)=0\},$$
\begin{multline*}
J_\Om^2(z;X,Y)=\sup\Big{\{}\big{|}\sum_{j,k=1}^n\frac{\partial^2f}{\partial\zeta_j\partial\zeta_k}(z)X_jY_k\big{|}:\\
f\in L^2_h(\Om),||f||_D\leq 1,f(z)=0,f'(z)=0\Big{\}}.
\end{multline*}

It is well-known (see e.g. \cite{JP,KY,Yoo,Zha}) that
$J_\Om^0$ is equal to the Bergman kernel $K_\Om$ on the diagonal,
and $$\ds\beta^2_{\Om}(z;X)=\frac{J^1_\Om(z;X)}{J^0_\Om(z)},\quad
b_\Om(z;X,Y)=2-\frac{J^0_\Om(z)J_\Om^2(z;X,Y)}{J_\Om^1(z;X)J_\Om^1(z;Y)},$$
where $b_\Om$ is the bisectional curvature of $\beta_\Om$
(we assume that $K_D(z)\neq 0$ and $\beta_\Om(z;Z)\neq 0$ if $Z\neq 0$).
Recall that $s_\Om(z;X)=b_\Om(z;X,X)$ and
$$r_\Om(z;X)=\sum_{j=1}^n b_\Om(p;e_j,X),\quad t_\Om(z)=\sum_{j=1}^n r_\Om(p;e_j),$$
where $e_1,\ldots,e_n$ is any basis of $\C^n$ ($r_\Om$ and $c_\Om$
do not depend on the choice of the basis).

\begin{thm}\label{curv} If $(\ast)$ holds, $\var_k\neq 1,$ and $\Om$ is pseudoconvex,
then there exist a constant $C>0$ and a neighborhood $U$ of $p$ such that
$$|s_\Om(z;X)+4/(n+1)|,|r_\Om(z;X)+1|,|t_\Om(z)+n|
\le C\delta_\Om(z)^{\var_k}$$
for $z\in\Om\cap U,\,X\in(\C^n)_\ast.$ Moreover, the estimates from Theorem \ref{nik}
remain true for $\frac{\beta_\Om}{\sqrt{n+1}}$ instead of $\kappa_{\Om}.$
\end{thm}

The above theorem is still valid in the $\CC^2$-smooth case, replacing the term
$C\delta_\Om(z)^{\var_k}$ by any positive number (for the bounded case see also \cite{KY} and the references therein).

When $\Om$ is a $\CC^{3}$- or $\CC^{4}$-smooth strictly pseudoconvex domain,
Theorem \ref{curv} is contained in \cite{Zha}.

To prove Theorem 4, we need the following localization result:

\begin{thm}\label{loc} Let $p$ be a strictly pseudoconvex boundary
point of a pseudoconvex domain $\Om$ in $\C^n$, and let
$V$ be a bounded neighborhood of $q.$ There exist a neighborhood $U\subset
V$ of $p$ and a constant $C>0$ such that if $z\in\Om\cap U$ and $X,Y\in(\C^n)_\ast,$ then
$$1\ge\frac{J^0_\Om(z)}{J_{\Om\cap V}(z)},\frac{J^1_\Om(z;X)}{J^1_{\Om\cap V}(z;X)},
\frac{J^2_\Om(z;X,Y)}{J^2_{\Om\cap V}(z;X,Y)}\ge 1+C\delta_\Om(z)\log\delta_\Om(z).$$
\end{thm}

Theorem \ref{loc} is proved in \cite{Nik} for $J^0$ and $J^1.$ The proof for $J^2$ is similar
and we skip it. The proof in \cite{Nik} shows that Theorem \ref{loc} remains
true under the following general hypothesis: there exist neighborhoods
$V_3\subset V_2\Subset V_1\subset V$ of $p$ and a constant $C>0$ such that for any
$q\in\partial D\cap V_2$ there is a peak function $\theta_q$ for
$D\cap V_1$ at $q$ with
$$1-|\theta_q(z)|\le C||z-q||,\ z\in D\cap V_3,$$
$$\sup\{|\theta_q(z)|:q\in\partial D\cap V_2,z\in D\cap V_1\setminus V_2\}<1.$$
For example, this holds if $p$ is a locally convexifiable point of finite ty\-pe $2m$
(see \cite{DFW}); then $1-|\theta_q(z)|\ge c||z-q||^{2m}$ with $c>0$ independent of $q.$

We will also use the following estimates, say (E), (see \cite{Zha}):
$$2-2\frac{n+3}{n+1}\sigma^{-4n}_\Om(z)\le s_\Om(z;X)\le2-2\frac{n+3}{n+1}\sigma^{4n}_\Om(z),\footnote
{In \cite{Zha}, the number $n+3$ is replaced by $n+2$ because of the assumption there
that the sectional curvature of the unit ball in $\mathbb{C}^n$ is $-\frac{2}{n+1},$ not $-\frac{4}{n+1}.$}$$
$$(n+1)-(n+2)\sigma^{-2n}_\Om(z)\le r_\Om(z;X)\le(n+1)-(n+2)\sigma^{2n}_\Om(z).$$
Therefore,
$$n(n+1)-n(n+2)\sigma^{-2n}_\Om(z)\le t_\Om(z)\le n(n+1)-n(n+2)\sigma^{2n}_\Om(z).$$

\noindent{\it Proof of Theorem \ref{curv}.} We may choose a neighborhood $V$ of $p$ so that $\Om\cap V$
is a strictly pseudoconvex domain. Since $V$ can be chosen bounded, Theorem \ref{loc} provides a
neighborhood $U'\subset V$ of $p$ and a constant $C'>0$ such that if $z\in\Om\cap U'$ and $X,Y\in(\C^n)_\ast,$ then
$$|b_\Om(z;X,Y)-b_{\Om\cap V}(z;X,Y)|,\ \left|\frac{\beta_\Om(z;X)}{\beta _{\Om\cap V}(z;X)}-1\right|
\le-C'\delta_\Om(z)\log \delta_\Om(z).$$
It remains to apply the estimates (E) and Corollary \ref{cor}
for $\Omega\cap V$ instead of $\Omega,$ and then Theorem \ref{main}.\quad\qed
\smallskip

The above proof implies that Theorem \ref{curv} remains true in the $\CC^{3,1}$-smooth case ($\var_k=1$)
with $\delta_\Om\log\delta_\Om$ instead of $\delta_\Om.$ It is natural to ask whether the log term can be canceled.
The answer would be affirmative if Theorem \ref{loc} holds with no log term.

Note also that the proof of Theorem \ref{curv} can be adapted to other Bergman invariants which can be expressed by
extremal functions. For example, let
$$J_\Om(z)=\frac{\Delta(z)}{K_\Om(z)}$$
be the canonical Bergman invariant, where $\Delta_\Om(z)$ is the determinant of the quadratic form $\beta_\Om(z;\cdot)$
(see e.g. \cite{Eng,KY}). Then, under the hypo\-theses of Theorem \ref{curv}, for $z$ near $p$ we have that
$$|n!J_\Om(z)-(n+1)\pi^n|=O(\delta_\Om(z)^{\var_k}).$$
$\CC^2$-smooth case for bounded $\Omega$ is treated in \cite{KY}; then $o(1)$ appears.
\smallskip

\noindent{\bf Acknowledgement.} The authors would like to thank the referee for many helpful remarks that essentially improved
the presentation of the paper.

\end{document}